
\baselineskip=14pt
\parskip=10pt

\magnification=\magstephalf
\def\M{{\cal M}}

\def\1{{\overline{1}}}
\def\2{{\overline{2}}}
\parindent=0pt
\overfullrule=0in

\def\frac#1#2{{#1 \over #2}}

\bf
\centerline
{
A Motivated Rendition of the Ellenberg-Gijswijt Gorgeous proof that the Largest Subset of $F_3^n$ with No
}
\smallskip
\centerline
{
Three-Term Arithmetic Progression is $O(c^n)$, with $c=\root 3 \of {(5589+891\,\sqrt {33})}/8=2.75510461302363300022127...$
}
\rm
\bigskip
\centerline
{\it By  Doron ZEILBERGER}
\bigskip
Let $F_3:=\{0,1,2\}$ be the field of integers modulo $3$, and let ${{n} \choose {k}}_2$ be the {\it trinomial coefficient}, defined as
the coefficient of $x^k$ in $(1+x+x^2)^n$. As usual, the number of elements of a finite set $S$ will be
denoted by $|S|$.

Inspired by the Croot-Lev-Pach [CLP] breakthrough, Jordan Ellenberg and Dion Gijswijt[EG] have recently amazed the combinatorial world by proving

{\bf Theorem.} ([EG]) Let $A$ be a subset of $F_3^n$ such that the equation
$$
a+b+c=0 \quad, \quad (a,b,c \in A)
$$
has no solutions except the trivial $a=b=c$. Then
$$
|A| \, \leq \, 3 \sum_{k=0}^{\lfloor \frac{2n}{3} \rfloor}  {{n} \choose {k}}_2 \quad.
$$
They then went on to show (using more-advanced-than-necessary probability theory [``large deviations'']) that
$|A|= O(2.75510461302363300022127...^n)$, but as observed by Terry Tao ([T]), this can be derived in a more
elementary way, only using Stirling's approximation of $n!$ and the (very simple) discrete Laplace method, as outlined, for example,
by Knuth in [K] (pp. 65-67).

The reason that their result was such a sensation was that many smart people tried very hard to improve the $o(3^n)$ result
proved in 1982, by Tom Brown and Joe Buhler, that was improved, in 1995, to $(3^n/n)$, by Roy Meshulam, and the 
current record (before [EG])
was $O(3^n/n^{1+\epsilon})$ by Michael Bateman and Netz Hawk Katz that was considered ``significant'' enough to
be accepted by the ``prestigious'' {\it Journal of the American Mathematical Society}.
(See [EG] for references).

This is  reminiscent of the long-standing challenge to improve $|\sum_{i=1}^{n} \mu(i)|=O(n^{1-\epsilon})$ to
$|\sum_{i=1}^{n} \mu(i)|=O(n^c)$, for {\it some}  $c<1$ (even $c=1-10^{-100000}$). $c=\frac{1}{2}+\epsilon$ would get you
a million dollars, but any $c<1$ would be a major breakthrough.

The article [EG], while better-written than \%99 of mathematical papers, is still suboptimal, since
it suffers from mathematicians' bad habit to hide their motivation. The account below
is just a motivated, {\it top down}, rendition of the beautiful [EG] proof, 
aimed at the {\it proverbial} {\bf smart freshman} (who took basic linear algebra).

\vfill\eject

{\bf Motivated Proof}

We need an {\it upper bound} for $|A|$. Since the {\it polynomial method} 
and {\it linear algebra} are such powerful tools, let's try to find some {\it vector space} of polynomials
whose dimension can be bounded from {\it below} by some expression involving $|A|$, and of course $n$, and
possibly another natural parameter, $d$, that at the end of the day can be chosen optimally in terms of $n$.

What can be more natural than the vector space of polynomials in $n$ variables, $x_1, \dots, x_n$ on $F_3^n$?
This vector space has a natural basis consisting of the $3^n$ {\it monomials}
$$
\{ x_1^{\alpha_1} \cdots  x_n^{\alpha_n}  \quad | \quad 0 \leq \alpha_i \leq 2 \} \quad ,
$$
and hence the dimension of this space is $3^n$. Also natural are the subspaces $\M(n,d)$ of polynomials
of (total) degree $\leq d$, whose natural basis  is the set of monomials
$$
M(n,d):= \{ x_1^{\alpha_1} \cdots  x_n^{\alpha_n}  \, | \quad 0 \leq \alpha_i \leq 2 \quad, \quad \alpha_1 + \dots + \alpha_n \leq d \} \quad ,
$$
and hence the dimension of $\M(n,d)$, alias $|M(n,d)|$, is given explicitly by
$\sum_{i=0}^{d} {{n} \choose {i}}_2$.

[When you expand $(1+x+x^2)^n$ you have to decide for each factor $(1+x+x^2)$
whether it is $x^0$, $x^1$, or $x^2$, giving a term $x^{\alpha_1 + \dots + \alpha_n}$, and the number of 
such tuples $(\alpha_1, \dots, \alpha_n)$
that add-up to $i$ is the coefficient of $x^i$ in $(1+x+x^2)^n$, that is ${{n} \choose {i}}_2$].

Since we want to find a vector space, $V$, whose dimension can be bounded in terms of $|A|$ (our object of desire),
the {\it first try} would be to consider the subspace of $\M(n,d)$  of polynomials {\it vanishing} on $A$, that would entail
$\dim V \, \geq \, |M(n,d)| - |A|$, leading to $|A| \geq |M(n,d)|- \dim V$. Alas, this is a {\it lower} bound for $|A|$, while we are
after an {\it upper} bound.

So the next thing to try (and it works!) is to consider the subspace of $\M(n,d)$ of polynomials that vanish
on the {\bf complement} of $A$, $F_3^n \backslash A$
$$
V:= \{P(x_1, \dots, x_n) \, | \, degree(P) \leq d \quad, \quad P(x)=0 \quad \hbox{ for all $x \in F_3^n \backslash A$} \} \quad ,
$$
entailing the following bound 
$$
\dim V \geq |M(n,d)| -| F_3^n \backslash A|=|M(n,d)| -(3^n- |A|)=|A| - (3^n -|M(n,d)|) \quad,
$$
that gives the {\bf upper bound} $|A| \leq \dim V + (3^n - |M(n,d)|)$.

How can we bound $\dim V$? All the members of $V$ are polynomials that vanish on $F_3^n \backslash A$, hence their
{\bf supports} are all subsets of $A$. If $P \in V$ has a support of maximal size, let's call it $\Sigma$, then
$|\Sigma| \geq \, dim \, V$.
Indeed, suppose that $|\Sigma|<\, \dim\, V$. Then there would be a {\bf non-zero} member $Q \in V$ that vanishes on $\Sigma$.
Since $Q$ is not identically zero, there is a point outside $\Sigma$ in which $Q$ is non-zero, while
$P$ must be $0$ (since it is $0$ outside its support). Hence $P+Q$ is non-zero on $\Sigma$ and that extra
point, and hence its support is strictly larger than $\Sigma$ contradicting the assumption that $P$ was
a member with maximal support.

So we have the bound

$$
|A| \, \leq \, |\Sigma| + (3^n - |M(n,d)|) \quad .
$$
It remains to say something about the maximal size of the supports of members $P \in V$. 

So far this is true for {\it any} subset $A \subset F_3^n$. It is time to take advantage of the fact that
it can never happen that $a_0+b_0+c_0$=0, with $a_0,b_0,c_0 \in A$ and $b_0 \neq c_0$. 

Let's define a set $S$ by
$$
S:=\{-b_0-c_0 \, | \, b_0,c_0 \in A \quad , \quad b_0 \neq c_0 \, \} \quad .
$$

$S$ is {\bf disjoint} from $A$, hence is a {\bf subset}  of $F_3^n \backslash A$. So we know that
{\it every} $P \in V$ vanishes in $S$, i.e.
$$
P(-b_0-c_0)=0, \quad \hbox{whenever $b_0, c_0 \in A$ and $b_0 \neq c_0$} \quad .
\eqno(ZeroCondition)
$$
Consider the $|A|$ by $|A|$ matrix whose rows and columns are indexed by the members of $A$, and
whose $(b_0,c_0)$ entry is $P(-b_0-c_0)$. By Eq. $(ZeroCondition)$, this is a {\bf diagonal matrix}.

On the other hand, the  polynomial $P(-b-c)$, viewed as a polynomial of {\bf total degree} $\leq d$  in the $2n$ variables
$b_1, \dots, b_n; c_1, \dots,c_n$, is a sum of
{\bf monomials } of the form
$$
(b_1^{\beta_1} \cdots b_n^{\beta_n})\cdot (c_1^{\gamma_1} \cdots c_n^{\gamma_n}) \quad,
$$
where $\beta_1 + \dots + \beta_n +  \gamma_1 + \dots + \gamma_n \leq d$.

Each and every such monomial can be written either as $m(b)m'(c)$ with $\deg m(b) \leq d/2$
or $m(c) m'(b)$ with $\deg m(c) \leq d/2$
[If  $n$ married (heterosexual) couples are given $\leq d$ ice-creams either the men have $\leq d/2$ of them or
the women (or both, in which case you can split them)].

Collecting terms, we get the {\bf crucial observation} (due to [CLP]) that, for every polynomial $P$, of degree $\leq d$,
there exist polynomials $F_m$ (one for each monomial $m$ of degree $\leq d/2$) such that
we can write
$$
P(-b-c)= \sum_{m \in M(n,d/2)} m(b) F_m(c) \, + \, \sum_{m \in M(n,d/2)} m(c) F_m(b) \quad.
\eqno(CLP)
$$

{\bf Plugging-in} $b=b_0,c=c_0$ into Eq. $(CLP)$ yields
$$
P(-b_0-c_0)= \sum_{m \in M(n,d/2)} m(b_0) F_m(c_0) \, + \, \sum_{m \in M(n,d/2)} m(c_0) F_m(b_0) \quad.
$$
Hence our diagonal matrix (whose $(b_0,c_0)$-entry is $P(-b_0-c_0)$) is a sum of
$2 |M(n,d/2)|$ matrices (two for each monomial $m \in M(n,d/2)$).
The summand, the matrix whose $(b_0,c_0)$ entry is $m(b_0) F_m(c_0)$
has rank 1 (since all rows (and all columns) are proportional to each other). Ditto for $m(c_0) F_m(b_0)$.
Hence that diagonal matrix is a sum of $2\,|M(n,d/2)|$ rank-one matrices, and hence its rank is $\leq 2 |M(n,d/2)|$.
Hence that matrix can have at most $2 |M(n,d/2)|$ non-zero diagonal entries, and hence
$P(-b_0-b_0)=P(b_0)$ is non-zero for {\bf at most} $2 |M(n,d/2)|$ members of $b_0 \in A$, and hence the size of the
{\bf support} of every $P \in V$ is  at most $2 |M(n,d/2)|$. In particular $|\Sigma| \leq 2 |M(n,d/2)|$.

We now got a family of {\bf explicit} upper bounds
$$
|A| \leq 2\,|M(n,d/2)|  + 3^n - |M(n,d)| \quad ,
$$
valid for {\it every} $d$. It turns out (and is easy to check on the computer) that taking
$d=\frac{4}{3}n$ will make it as small as possible. For the sake  of convenience let's assume that
$n$ is a multiple of $3$. We get
$$
|A| \leq 2\,|M(n,\frac{2}{3} n)|  + 3^n - |M(n, \frac{4}{3}n)| \quad .
$$
Since, by symmetry  ${{n} \choose {2n-k}}_2= {{n} \choose {k}}_2$, we have:
$$
3^n- |M(n,\frac{4}{3}n)|= 3^n- \sum_{k=0}^{\frac{4}{3}n}  {{n} \choose {k}}_2 =
\sum_{k=0}^{2n}  {{n} \choose {k}}_2 - \sum_{k=0}^{\frac{4}{3}n}  {{n} \choose {k}}_2 =
\sum_{k=\frac{4}{3}n+1}^{2n}  {{n} \choose {k}}_2
$$
$$
= \, \sum_{k=0}^{\frac{2}{3}n - 1}  {{n} \choose {k}}_2 \, = \,
\sum_{k=0}^{\frac{2}{3}n}  {{n} \choose {k}}_2 - {{n} \choose {\frac{2}{3}n}}_2 \quad .
$$
Hence
$$
|A| \leq 3\, \sum_{i=0}^{\frac{2}{3}n} {{n} \choose {i}}_2 -  {{n} \choose {\frac{2}{3}n}}_2 
\leq \, 3\, \sum_{i=0}^{\frac{2}{3}n} {{n} \choose {i}}_2 \quad .
$$
It is easy to see that this is $\leq C  {{n} \choose {\frac{2}{3}n}}_2$ for some positive constant $C$, so it remains
to find the asymptotics of ${{n} \choose {\frac{2}{3}n}}_2$.

{\bf  Asymptotics}

[EG] used the sledge-hammer of ``large deviations'', but
as noticed in [T], the asymptotics can be derived by {\bf purely elementary methods}.
An even better ({\it and even more elementary}!) way to find the asymptotics is to use the {\bf Almkvist-Zeilberger Algorithm} [AZ],
as implemented in the Maple package  

{\tt http://www.math.rutgers.edu/\~{}zeilberg/tokhniot/EKHAD} \quad .

Since  ${{3n} \choose {2n}}_2$ is the constant term of $(1+x+x^2)^{3n}/x^{2n}$, typing in {\tt EKHAD}

{\tt AZd((1+x+x**2)**(3*n)/x**(2*n+1),x,n,N)[1];}

immediately yields the linear recurrence operator annihilating the sequence $d(n):={{3n} \choose {2n}}_2$,
viz. that $d(n)$ satisfies the second order linear recurrence equation with polynomial coefficients
\vfill\eject
$$
243\, \left( 3\,n+5 \right)  \left( 3\,n+2 \right)  \left( 11\,n+20 \right)  \left( 3\,n+4 \right)  \left( 1+3\,n \right)  \left( n+1 \right) 
d \left( n \right) 
$$
$$
-18\, \left( 3\,n+5 \right)  \left( 1+2\,n \right)  \left( 3\,n+4 \right)  \left( 759\,{n}^{3}+2898\,{n}^{2}+3505\,n+1350 \right) d \left( n+1
 \right) 
$$
$$
+16\, \left( 5+4\,n \right)  \left( 3+2\,n \right)  \left( 1+2\,n \right)  \left( 11\,n+9 \right)  \left( 7+4\,n \right)  \left( n+2 \right) d
 \left( n+2 \right) \, = \, 0 \quad .
$$
By the Poincar\'e lemma, $d(n)$ is asymptotic (ignoring $n^{\alpha}$ terms), (taking the leading coefficient in $n$, namely  $n^6$, in the above
recurrence), to the solution, $d_0(n)$ of the
linear recurrence with {\bf constant coefficients}
$$
19683\,d_{{0}} \left( n \right) -22356\,d_{{0}} \left( n+1 \right) +1024\,d_{{0}} \left( n+2 \right) \, = \, 0 \quad ,
$$
whose largest characteristic root is the root of 
$$
1024\,{N}^{2} -22356\,N +19683=0 \quad ,
$$
that happens to be ${\frac {5589}{512}}+{\frac {891}{512}}\,\sqrt {33}= 20.912901011846452219 \dots $, and taking the cubic root, we get
that $|A|=O(\alpha^n)$ where
$$
\alpha \, = \,  2.7551046130236330002\dots \quad .
$$
Using the Maple package {\tt http://www.math.rutgers.edu/\~{}zeilberg/tokhniot/AsyRec.txt}, one can get the more precise
asymptotics $|A| \leq C \alpha^n \frac{1}{\sqrt{n}}$, for some $C$. In fact, we have:
$$
|A| \, \leq \,  
 3.3267627467425979588\,\cdot { (2.7551046130236330002\dots)}^{n} \cdot {\frac {1}{\sqrt {n}}} \cdot
$$
$$
\left(  1 \, -\,  5.1543714155636062458\,{n}^{-1}+ 90.161538946865747706\,{n}^{-
2}- 2646.8299396834595447\,{n}^{-3} + O(n^{-4}) \right) 
\quad .
$$

{\bf ${\bf  F_q^n}$}

As pointed out in [EG] analogous arguments can be applied for $F_q^n$ for {\it any} prime power $q$.
The same elementary argument described in [T] works in general (yielding the same answers given by
large deviations), and our approach, via the Almkvist-Zeilberger algorithm, works also works well.

The [EG] upper bound for general $q$ is expressible as the coefficient of $z^{(q-1)n/3}$ in the rational function
$$
(1+z+ \dots + z^{q-1})^n \cdot \frac{2+z}{1-z} \quad .
$$
For every given $q$, the Almkvist-Zeilberger algorithm produces a rerurrence (if $q-1$ is not divisible by $3$ one
has to replace $n$ by $3n$), from which the asymptotics can be deduced as above. Alternatively, one can express
that quantity as a contour-integral and use Laplace's method for integrals. The advantage of the latter method is
that one can handle {\it all} $q$ in one stroke, i.e. leave $q$ {\bf symbolic}.

For the record, here are the growth constants for primes and prime powers $ 4 \leq q \leq 31$.

$q=4 \, : \, 3.610718613276039349\dots $ \quad; \quad 

$q=5 \, : \, 4.461577765702577811\dots$ \quad; \quad 

$q=7 \, : \, 6.156204863216738416 \dots$ \quad; 

$q=8 \, : \, 7.0015547549940074584$ \quad;\quad

$q=9 \, :  \, 7.846120582585805712\dots$ \quad;

$q=11 \, : \, 9.533685392075550992\dots$ \quad;

$q=13 \, :  \, 11.21990798911487743\dots$ \quad;

$q=16 \, : \,  13.74776213458745700\dots$ \quad ;

$q=17 \, : \,  14.590117162\dots   $ \quad ;

$q=19 \, : \,  16.274551068400264\dots $ \quad ;

$q=23 \, : \,  19.6426364587288\dots $ \quad ;

$q=25 \, : \, 21.3264083101\dots  $ \quad ;

$q=27 \, : \,   23.010051182485787\dots    \dots  $ \quad ;

$q=29 \, : \,  24.69359086763659\dots \dots  $ \quad ;

$q=31 \, : \,  26.3770467097314914 \dots  $ \quad .

{\bf References}

[AZ] Gert Almkvist and Doron Zeilberger, {\it The Method of Differentiating Under The Integral Sign},
J. Symbolic Computation {\bf 10}(1990), 571-591. Available on-line from \hfill\break
{\tt http://www.math.rutgers.edu/\~{}zeilberg/mamarim/mamarimhtml/duis.html} \quad .

[CLP] Ernie Croot, Vsevolod Lev, Peter Pach, {\it Progression-free sets in $Z_4^n$ are exponentially small},
May 5, 2016, {\tt http://arxiv.org/abs/1605.01506} \quad .

[EG] Jordan S. Ellenberg and Dion Gijswijt, {\it On large subsets of $F_q^n$ with no three-term arithmetic progression},
May 30, 2016, {\tt https://arxiv.org/abs/1605.09223} \quad . 

[K] Donald E. Knuth, {\it ``The Art of Computer Programming, vol. III: Sorting and Searching''}, Addison-Wesley, 1973.

[T] Terence Tao, {\it A symmetric formulation of the Croot-Lev-Pach-Ellenberg-Gijswijt capset bound},
May 18, 2016, {\tt https://terrytao.wordpress.com/tag/polynomial-method/} \quad .     

\vfill\eject

\bigskip
\bigskip
\hrule
\bigskip
Doron Zeilberger, Department of Mathematics, Rutgers University (New Brunswick), Hill Center-Busch Campus, 110 Frelinghuysen
Rd., Piscataway, NJ 08854-8019, USA. \hfill \break
zeilberg at math dot rutgers dot edu \quad ;  \quad {\tt http://www.math.rutgers.edu/\~{}zeilberg/} \quad .
\bigskip
\hrule

\bigskip
Exclusively published in The Personal Journal of Shalosh B. Ekhad and Doron Zeilberger  \hfill \break
({ \tt http://www.math.rutgers.edu/\~{}zeilberg/pj.html})
and {\tt arxiv.org} \quad . 
\bigskip
\hrule
\bigskip
{\bf  Written:  July 6, 2016}

\end